\documentclass[11pt]{amsart}
\usepackage{graphicx}
\usepackage{epsfig}
\usepackage[left=1.25in,top=1in,right=1.25in]{geometry}
\usepackage{setspace}

\theoremstyle{definition}

\newcommand{\ignore}[1]{}

\begin{document}

\title{N-Person Envy-Free Chore Division}

\author{Elisha Peterson}
\address{Department of Mathematics\\ United States Military Academy\\West Point, NY 10996 }
\email{elisha.peterson@usma.edu}

\author{Francis Edward Su}
\address{Department of Mathematics\\ Harvey Mudd College\\ Claremont, CA  91711}
\email{su@math.hmc.edu}



\thanks{Research partially supported by NSF Grant DMS-0701308 (Su).}

\dedicatory{\ \ \ }
\maketitle

\section{Introduction} 
\label{Intro}

In this paper we consider the problem of {\em chore division}, which is
closely related to a classical question, 
due to Steinhaus \cite{steinhaus}, of how to cut a cake fairly.  We focus
on {\em constructive} solutions, i.e., those obtained via a 
well-defined procedure or algorithm.
Among the many notions of fairness 
is {\em envy-freeness}: an {envy-free} 
cake division is a set of cuts and an allocation of the pieces
that gives each person what she feels is the {\em largest}
piece.  Much progress has been made on finding constructive algorithms 
for achieving envy-free cake divisions; a landmark result was that of 
Brams and Taylor \cite{bramstaylormonthly}, who gave the first general $n$-person procedure.

In contrast to cakes, which are desirable, the dual problem of chore division 
is concerned with dividing an object deemed {\em undesirable}.
Here, each player would like to receive what he considers to be
the {\em smallest} piece, of say, a set of chores.
This problem appears to have been first introduced by 
Martin Gardner in \cite{gardner}.  
Oskui (see \cite{robertsonwebb}) referred to it as the {\em dirty work
problem} and gave the first discrete and moving-knife solutions for
exact envy-free chore division among 3 people.  Peterson and Su \cite{petersonsu}
gave the first explicit 4-person moving-knife procedure for chore division, adapting ideas of 
Brams, Taylor, and Zwicker \cite{bramstaylorzwicker} for cake-cutting.

The purpose of this article is to give a general $n$-person solution to the chore division problem.
Su \cite{su} gives an $n$-person chore division algorithm but it only yields an $\epsilon$-approximate solution after a finite number of steps.
Brams and Taylor suggest in \cite{bramstaylor} how cake-cutting methods
could be adapted to chore division without working out the details, 
and our algorithm owes a great debt to their ideas.
But we also show where some new ideas are needed, and why the chore division problem is not exactly a dual or straightforward extension of the cake-cutting problem.

\section{Some Key Ideas}

We assume throughout this paper that chores are infinitely divisible. 
This is not unreasonable as a
finite set of chores can be partitioned by 
dividing up each chore (e.g., a lawn to be mowed could
be divided just as if it were a cake), or dividing the time spent on them. 
For ease of expression, we shall call
the set to be divided a {\em cake}, rather than a {\em set of chores}.
Implicit in this is the assumption that the players
desire the {\em smallest, rather than largest, piece of cake}.

We remind the reader that players may have different preferences over pieces of cake (indeed that is what makes the problem interesting).  More formally, each player $i$ has a 
{\em measure} $\mu_i$ that describes what value $\mu_i(A)$ the player assigns to a piece of cake $A$.  (The cake pieces in our construction will always be measurable).  Such measures are {\em additive}, meaning that no value is created or destroyed by cutting or lumping pieces together.

Before addressing the chore division problem, we wish to highlight a couple of key ideas from 
$n$-person cake-cutting algorithm of Brams and Taylor, and discuss the analogous ideas in the chore division context.

\begin{itemize}
\item {\bf Trimming to Create Ties.}
Given a division of cake into several pieces, a player $B$ can create a ``tie'' in what she considers the largest piece by trimming the largest piece (the trimmings are temporarily set aside).  This way, if another person $C$ chooses a piece before her, player $B$ will still have one of her top two choices for largest piece available to choose, so she will not envy $C$ who chose before her.

\item {\bf Irrevocable Advantage.}  
Suppose that the cake has been allocated, except for trimmings, in such a way 
that player $A$ receives a piece that he believes to be $\epsilon$ bigger in his measure 
than the piece $B$ received.  Suppose also that the trimmings are, in $A$'s estimation, of size less than $\epsilon$.  Then no matter how much of the trimmings are given to $B$, player $A$ will never envy her.  We say that $A$ has an {\em irrevocable advantage} over $B$.
\end{itemize}

For chore division, the idea corresponding to trimming 
is the idea of ``adding back'' or {\em augmentation} from a set of {\em reserves}.  We will also use an analogous concept of an irrevocable advantage.  Both ideas are present in the $3$-person chore division algorithm of Oskui (see \cite{robertsonwebb}), and we will use them repeatedly in our $n$-person algorithm.

\begin{itemize}
\item {\bf Augmentation with Reserves.}
Suppose a player has set aside, in advance, a set of cake we will call her {\em reserves}.
Given a division of cake, a player $B$ can
create a tie for {\em smallest} piece (in her opinion) 
by augmenting the smallest piece with cake from her reserves.
This is to ensure that if another player $C$ chooses a piece before her, she will still have one of her two smallest pieces available to choose from, so that $B$ will not envy $C$.

\item {\bf Irrevocable Advantage.}
Suppose that the cake has been allocated, except for unused reserves, in such a way 
that player $A$ receives a piece that he believes to be $\epsilon$ {\em smaller} in his measure than  the piece that $B$ received.  Suppose that the unused reserves are, in $A$'s estimation, of size less than $\epsilon$.  Then no matter how much of the unused reserves are given to $A$, player $A$ will not envy player $B$.  We say that $A$ has an {\em irrevocable advantage} over $B$.
\end{itemize}

Of course, some issues that we will have to resolve in our algorithm are:
(1) how to create enough reserves for players to use and (2) how to deal with unused reserves.
Irrevocable advantages can be used for (2) when the reserves are small enough.

\section{An $N$-Person Envy-Free Chore Division Procedure} 
\label{NPerson}

We now construct our $n$-person chore-division procedure. 
Unlike the 4-person moving-knife scheme of Peterson and Su
\cite{petersonsu}, ours is a discrete procedure (involving no continuous evaluations of pieces like moving-knife methods require).  
And while their $4$-person procedure 
possessed a natural set of reserves due to the initial trimming, for our $n$-person procedure 
we need to carefully create enough reserves for use by specific players.



A brief sketch of our procedure runs as follows.  Let one player divide
the cake and allocate the pieces.   As long as there are objections, 
we shall iterate a procedure that gives an envy-free allocation of part
of the cake (Steps 1-9), and also gives a player
who objected an irrevocable advantage over another player with respect
to the part of the cake that has not yet been allocated (Steps 10-15).  With enough
iterations on the leftovers, there will be enough players with irrevocable advantages 
to allow the allocation of the remainder of the cake (Steps 16-20).

As we noted, our method closely mirrors Brams and Taylor's $n$-person cake-cutting
procedure \cite{bramstaylormonthly}, 
but differs from theirs in using augmentation (rather than trimming) and the creation of reserves.
For ease of comparison with Brams and Taylor's cake-cutting procedure,
we include step numbers in our procedure that correspond their step
numbers \cite{bramstaylormonthly}.  The significant departures
occur in Steps 6.1, 7.1, 7.3, 8, and 10-15.  

Following their
example, we distinguish {\em rules} from {\em strategies} by placing
strategies in parentheses.  Again, for ease of expression, we refer to
the chore set as a cake, bearing in mind that each person wants
the piece he/she thinks is {\em smallest}.

We shall exhibit our procedure for $n=4$ case.  The generalization to more
players will be discussed subsequently.


\begin{quote}
\begin{enumerate}

\item[\bf Step 1.]
Let Player $2$ cut the cake into 4 pieces (that she considers equal), and 
then assign one piece to each player.

\item[\bf Step 2.]
Player 2 asks the other 3 players if anyone is envious.

\item[\bf Step 3.]
If no one has envy, then each keeps the piece he was given, and we are done.

\item[\bf Step 4.]
If someone has envy, say Player 1, let Player 1 choose two pieces (that
he thinks are not equal in size) and name them $A$ (for the larger piece)
and $B$ (for the smaller piece).

\end{enumerate}
\end{quote}
{\em Aside.}
The other pieces are reassembled for allocation later.  Note that 
Player 1 thinks $A$ is larger (hence less desirable) than $B$ but Player 2 thinks they are the
same size.

\begin{quote}
\begin{enumerate}

\item[\bf Step 5.]
Let Player $1$ name an integer $r\geq 11$ (chosen such that, 
even if $A$ were divided into $r$ pieces and the 8 smallest
pieces of $A$ were removed, he would still prefer $B$).

\end{enumerate}
\end{quote}

{\em Aside.}
This is possible because 
the union of the 8 smallest pieces is no larger than 8 times
the average size of all $r$ pieces.  Hence Player 1 can choose $r$ large
enough so that $8 \mu(A)/r < \mu(A) - \mu(B)$ where $\mu$ is Player 1's
measure.

\begin{figure}[htb]
  \begin{center}
    \scalebox{1}{\includegraphics{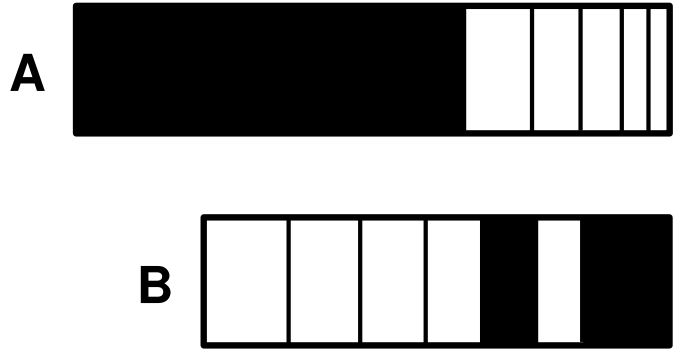}}
  \end{center}
  \caption{Steps 6-8:  Player 1's view of the cakes and piece sizes. Player 2 views all pieces as equal. 
$Z_1,Z_2,Z_3$ will either be the 3 shaded pieces in $A$ or the largest piece $A_1$ split into thirds.  
(The number of pieces $r$ was chosen large enough so that in one of these ways Player 1 would feel the $Z_i$ were strictly larger than the shaded pieces in $B$.)
If $A_1$ was split to form the $Z_i$, Player 2 may use the dotted pieces (her reserves) to equalize the $Z_i$ in her opinion.
In $B$, Player 3 chooses, say, the striped pieces as his reserves.  After these are chosen Player 1 chooses the smallest 3 remaining (here, the shaded pieces) in $B$ to be $Y_1,Y_2,Y_3$, then augments them from his reserves (unshaded pieces of $B$) to make them equal in his opinion.  
}
  \label{chore-pieces}
  \bigskip
\end{figure}

\begin{quote}

\begin{enumerate}

\item[\bf Step 6.]
Player $2$ divides each of $A$ and $B$ into $r$ sets (that
she considers equal).

\item[\bf Step 6.1.]  
From the pieces in $B$, 
let Player 3 choose two pieces (he considers largest).  This will be used
as Player 3's reserves.

\item[\bf Step 7.]  From the remaining pieces, Player $1$ chooses (what he
thinks are the smallest) $3$ sets in $B$, calling these $Y_1$, $Y_2$,
and $Y_3$.  

\item[\bf Step 7.1.]The rest of $B$ is set aside as Player 1's reserves.  If
necessary, Player $1$ uses his reserves to add to two of the $Y_i$ (to
make three of the $Y_i$'s equally sized--- the reserves will be enough because
they came from at least two pieces in $B$ that Player 1 feels is 
at least as large as each of the $Y_i$'s).

\item[\bf Step 7.2.]
If Player 1 considers the three largest pieces in $A$ all strictly 
larger than these pieces, the three largest pieces in $A$ 
are identified as $Z_1$, $Z_2$, and $Z_3$.

Otherwise, Player 1 cuts the largest piece in $A$ into three
(equal) pieces, calling them $Z_1$, $Z_2$, and $Z_3$ (which Player 1 will feel is
strictly bigger than each of the $Y_i$). 

\end{enumerate}
\end{quote}
{\em Aside.}
We show why one of the two cases in Step 7.2 must hold.
Let $\mu$ denote Player 1's measure, and sets $A_1,...,A_r$ and $B_1,...,B_r$ the
pieces of $A$ and $B$ arranged by decreasing $\mu$-size.
Since Player 3's reserves were chosen {\em before} the $Y_i$'s, the $Y_i$'s are among 
the 5 smallest sets $B_{r-4}, ... ,B_r$.
Suppose, to contradict both cases in Step 7.2, that both the following hold: 
(i) $\mu(B_{r-4}) \geq \mu(A_3)$ and
(ii) $\mu(B_{r-4}) \geq \mu(A_1)/3$.
From (i),
$$\mu(B_7 \cup \cdots \cup B_{r-4}) \geq \mu(A_3 \cup \cdots \cup A_{r-8}),$$
because there are
$r-10$ sets in each union and the 
smallest of the $B$ sets is at least as large as the largest of the $A$ sets.
From (ii),
$$\mu(B_1 \cup \cdots \cup B_6 ) \geq \mu(A_1 \cup A_2),$$ 
since each group of three $B$ sets is at least as large $A_1$ or $A_2$.
Taken together,
$\mu(B) \geq \mu(A_1 \cup \cdots \cup A_{r-8})$, contradicting our choice of $r$ in Step 5.

\begin{quote}
\begin{enumerate}
\item[\bf Step 7.3.]
The remaining pieces of $A$ are set aside for
Player 2's reserves.  If necessary,
Player 2 uses her reserves to add on to two of the $Z_i$ (to make all three $Z_i$'s equally sized--- this would only be needed in the second case of Step 7.2, and the 
reserves are enough because they came from
pieces in $A$ that were {\em each} the same size as all three $Z_i$'s combined).

\end{enumerate}
\end{quote}
{\em Aside.}
At this stage, Player $1$ believes $Y_1=Y_2=Y_3<Z_1,Z_2,Z_3$ (note
the strict inequality).  
Player $2$ believes $Z_1=Z_2=Z_3\leq Y_1,Y_2,Y_3$.  See Figure \ref{chore-pieces}.

\begin{quote}
\begin{enumerate}

\item[\bf Step 8.]
Player $3$ takes this collection of 6 pieces and, if necessary, augments
one of the pieces using his reserves 
(to make a two-way tie for {\em smallest} piece--- 
Player 3 has enough reserves because the 2nd-smallest piece cannot be larger than the 2nd-smallest $Y_i$, and Step 6.1 gave him reserves from two pieces of $B$, each of which he feels is larger than each $Y_i$ and larger than any possible augmentation of $Y_i$ by Player 1 in Step 7.1).

\item[\bf Step 9.]
Let the players choose in the order $4$-$3$-$2$-$1$, with Player 
$3$ required to take a piece he augmented if it is available.
Player $2$ must choose one of the $Z_i$'s, and Player $1$ must
choose one of the $Y_i$'s.

\end{enumerate}
\end{quote}
{\em Aside.}
This yields a partition $\{ X_1, X_2, X_3, X_4, L_1\}$ of the cake
such that the $X_i$'s are allocated in an envy-free fashion, and $L_1$ is
the leftover piece consisting of all cake not yet allocated.
Moreover, note that Player $1$ thinks his piece
is {\em strictly smaller} than Player $2$'s piece, say, 
by an amount $\epsilon$.

\begin{quote}
\begin{enumerate}

\item[\bf Step 10.]
Player $1$ names $s$ (such that $[\frac{1}{2}\mu(L_1)]^s<\epsilon$, 
where $\mu$ is Player 1's measure).  

\end{enumerate}
\end{quote}
{\em Aside.}
Steps 11-14 will result in Player 1 thinking that 
at least half of $L_1$ is allocated.  Hence $s$ represents the number 
of times to iterate Steps 11-14 to make the leftover piece smaller than 
$\epsilon$.

\begin{quote}
\begin{enumerate}

\item[\bf Step 11.]
Player $1$ cuts $L_1$ into $8$ (equal) pieces. Player $2$ sets aside (the
largest) 2 pieces for her reserves, 
and Player $3$ sets aside a piece from those
remaining (that he feels is largest of those remaining) for his reserves.  

\item[\bf Step 12.]
Of the remaining $5$ pieces, Player $2$ returns part of her
reserves to 2 of the pieces, if necessary (to create a 3-way tie
for the smallest piece).

\item[\bf Step 13.]
Player $3$ returns, if necessary, some from his reserves to one of the 5
pieces (to create a 2-way tie for the smallest piece--- his reserves will be sufficient
because his piece is at least as large as the 2nd-smallest piece).

\item[\bf Step 14.]
Let the players choose in the order $4$-$3$-$2$-$1$, with Players
$2$ and $3$ required to take a piece they augmented if one is
available.  
The chosen pieces are combined with the corresponding pieces from the players' earlier envy-free allocation 
to form a new envy-free partial allocation of cake $\{ X'_1, X'_2, X'_3, X'_4\}$, 
and the non-chosen pieces are lumped together to form a new smaller leftover piece $L_2$.

\item[\bf Step 15.]
Repeat steps 11-14 $(s-1)$ additional times, with each application of these
four steps applied to the leftover piece $L_2$ from the preceding application.

\end{enumerate}
\end{quote}

{\em Aside.}
According to Player $1$, the leftover 
$L_2$ is now smaller than $\epsilon$, so Player 1 feels that his portion together with $L_2$ (or any part of $L_2$) 
would still be smaller than Player 2's portion.  Thus, Player $1$ has an {\em irrevocable advantage} 
over player $2$ with respect to the leftovers.  We consider a subset of ordered pairs in 
$\{1,2,3,4\}\times\{1,2,3,4\}$ called $\mathcal{IA}$, to keep track of which player pairs have irrevocable advantages
with respect to the leftovers, and we place the ordered pair $(1,2)$ in $\mathcal{IA}$.   Note by definition, any further subdivision and allocation of pieces of $L_2$ with a smaller set of leftovers will not remove $(1,2)$ from $\mathcal{IA}$.

\begin{quote}
\begin{enumerate}

\item[\bf Step 16.]
Player $2$ cuts the remaining leftovers $L_2$ into 12 pieces (that she feels are the same size).

\item[\bf Step 17.]
Each player who agrees that all pieces are the same size is placed
in the set $A$.  Otherwise, players who disagree are placed in the set $D$.

\item[\bf Step 18.]
If $D\times A\subset \mathcal{IA}$, we divide the pieces among the
players in $A$, each receiving the same number of pieces, and we
are done.

\end{enumerate}
\end{quote}
{\em Aside.}
Players in $A$ do not envy each other since they agree each of 
the 12 pieces were the same size.  None of the players in $D$ envy
those in $A$ by the definition of $\mathcal{IA}$.  Players in $D$ do not
envy each other because they have not received any new pieces in this step.

\begin{quote}
\begin{enumerate}

\item[\bf Step 19.]
Otherwise, we choose the lexographically least pair $(i,j)$ from
$D\times A$ that is not in $\mathcal{IA}$, and return to step 4,
with Player $i$ in place of Player $1$, Player $j$ in place of
Player $2$, and $L_2$ in place of the cake.

\item[\bf Step 20.]
Repeat Steps 5-18.  

\end{enumerate}
\end{quote}
{\em Aside.}
Since each pass through Step 15 adds a new ordered pair to $\mathcal{IA}$ without removing any ordered pairs,
eventually $\mathcal{IA}$ will contain $D \times A$, no matter what $D \times A$ currently is.
When this occurs, the algorithm
concludes at Step 18 with an envy-free chore division of the entire cake.
This ends the procedure.

\bigskip
The extension of this procedure to $n$ players is tedious but straightforward, 
so we briefly mention the changes that ensue and leave
the verification to the reader.

Let $k=1 + 2 + \dots + (n-3)$.
The 8-piece criterion in Step 5 needs to be increased (to $k^2+4k+3$ pieces)
for later use in Step 7.
This necessitates an increase in $r$.
In Step 6.1, all players except 
Players 1, 2, and $n$ will choose, 
in reverse order, pieces of $B$ to form their reserves.  Specifically,
player $i$ will choose what she thinks are the $2(n-i)$ largest 
pieces of what remains in $B$; this ensures that these players have enough reserves for Step 8.
Thus the total number of pieces chosen for reserves in Step 6.1 is $2k$.

In Step 7, there need to be more $Y_i$'s and $Z_i$'s chosen ($k+2$ of each), 
so that after the augmentation in Step 8 
there are at least two $Y_i$'s and two $Z_i$'s untouched.  
Then Step 9 will allow players to choose in reverse order and still have one of the
$Z_i$'s available for player $2$ and one of the $Y_i$'s available for Player 1. 
Also, in Step 7.2, if Player 1 feels the biggest $k+2$ pieces of $A$ are not all larger
than the $Y_i$'s just chosen, then $r$ needs to be large enough 
($r\geq k^2+5k+5$) so that 
the largest piece of $A$ can be split evenly into $k+2$ pieces
all larger than the $Y_i$'s.  (Some care needs to be exercised here, since
the reserves are chosen {\em before} the $Y_i$'s, so Player 1
may only conclude that the $Y_i$'s are among the $2k+1$ smallest pieces of $B$.)
The size of $r$ also guarantees that Player 1 and 2 have enough reserves from
the remainders of $B$ and $A$.

In Step 8, there are enough reserves because each person feels his
reserves are enough to cover the correct number of pieces of $Y_i$ and any possible augmentations by Player $1$.
In Step 9, the players choose pieces in reverse order and players are required
to take a piece they augmented if available.  Steps 10-15 can be modified
analogously (let Player 1 start by cutting the leftovers into
$n^2-3n+4$ equal pieces) to accomodate more people 
in the iterative part of the procedure.

\section{Remarks}
We have shown explicitly 
how cake-cutting algorithms can be translated,
with complications, into exact envy-free chore-division algorithms.

There are two important features of this translation.  First, the use
of augmentation from reserves is very important, and we showed that the creation of such
reserves is generally possible.
Secondly, the notion of an irrevocable advantage for chores allows 
one to terminate what might otherwise be an infinite procedure.

Curiously, though our chore-division procedure
is more complicated than its cake-cutting
counterpart, it may converge faster and
require fewer cuts overall.  
For instance, for $n=4$, each pass through 
the iterative part of the procedure (Steps 10-15)
guarantees at least half the cake is apportioned
rather than only one-fifth.  The authors have noted that faster convergence
of chore division is also a feature of the $\epsilon$-approximate algorithm
of Su \cite{su}.

Note that unlike the 4-person moving-knife procedure of Peterson and Su \cite{petersonsu}
(which required at most 16 cuts), this $n$-person algorithm 
may take arbitrarily many cuts (and steps) to resolve, depending on player preferences, even for fixed $n$.  
So the number of steps is finite, but not bounded.
It remains an open question whether a {\em bounded} procedure
exists for either cake or chore division among $n$-people.

Given the numerous cuts needed to implement both the Brams-Taylor $n$-person 
cake-cutting procedure \cite{bramstaylor} 
and our $n$-person chore-division procedure, 
one may rightfully question their practicality.  There are two possible
responses.  First, the number of cuts can be reduced if one is 
willing to accept an $\epsilon$-approximate solution, by rotating players
through the roles of Steps 10-15, and quitting when satisfied.  Secondly, the 
construction of an initial solution, however complex, 
is always the first step towards finding useful simplifications.


\end{document}